\documentclass{article}
\usepackage{graphicx}
\usepackage{amsmath}
\usepackage{amssymb}
\usepackage{amsthm}

\newtheorem{thm}{Theorem}
\newtheorem{prop}{Proposition}
\newtheorem{cor}{Corollary}
\newtheorem{lem}{Lemma} 
\newtheorem{definition}{Definition}
\newtheorem{rem}{Remark}
\newtheorem{example}{Example}

\title{Bernstein Functions and Radial Limits of Prescribed Mean Curvature Surfaces}
\author{Mozhgan ``Nora'' Entekhabi \\ Department of Mathematics \\ Florida A \& M University \\ Tallahassee, FL 32307  \\ 
\& \\ 
Kirk E. Lancaster \\ Wichita, Kansas 67226}
\date{ }

\def\Real{{\rm I\hspace{-0.2em}R}}
\def\Natural{{\rm I\hspace{-0.2em}N}}

\newcommand\myeq{\mathrel{\overset{\makebox[0pt]{\mbox{\normalfont\tiny\sffamily def}}}{=}}}

\begin{document}

\maketitle

\vspace*{3mm}

\begin{abstract}
The radial limits at a point ${\bf y}$  of the boundary of the domain $\Omega\subset\Real^{2}$  of a bounded variational solution $f$  of Dirichlet or contact angle 
boundary value problems for a prescribed mean curvature equation are studied with an emphasis on the effects of assumptions about the curvatures of 
the boundary $\partial\Omega$  on each side of the point ${\bf y}.$  For example, at a nonconvex corner ${\bf y},$  we previously proved that all nontangential 
radial limits of $f$  at ${\bf y}$  exist; here we provide sufficient conditions for the tangential radial limits to exist, even when the Dirichlet data 
$\phi\in L^{\infty}(\partial\Omega)$  has no one-sided limits at ${\bf y}$  or the contact angle $\gamma\in L^{\infty}(\partial\Omega:[0,\pi])$  
is not bounded away from $0$  or $\pi.$  We also provide a complement to a 1976 Theorem by Leon Simon on least area surfaces. 
\end{abstract}

\section{Introduction}

Let $\Omega$  be a locally Lipschitz domain in $\Real^{2}$  and define $Nf = \nabla \cdot Tf = {\rm div}\left(Tf\right),$  where  $f\in C^{2}(\Omega)$  and 
$Tf= \frac{\nabla f}{\sqrt{1+\left|\nabla f\right|^{2}}}.$   Consider the Dirichlet problem
\begin{eqnarray}
\label{eq:D}
Nf & = & H(\cdot,f(\cdot))  \mbox{  \ in \ } \Omega    \\
f & = & \phi  \mbox{ \ on \ } \partial \Omega
 \label{bc:D}
\end{eqnarray}
and the contact angle problem
\begin{eqnarray}
\label{eq:C}
Nf & = & H(\cdot,f(\cdot))  \mbox{  \ in \ } \Omega    \\
Tf \cdot {\bf \nu} & = & \cos \gamma  \mbox{ \ on \ } \partial \Omega,
 \label{bc:C}
\end{eqnarray}
where $\phi:\partial\Omega\to\Real,$  $\gamma:\partial\Omega\to[0,\pi],$ and $H:\Omega\times\Real\to\Real$  are prescribed functions, 
$H({\bf x},t)$  is nondecreasing in $t$  for each ${\bf x}\in\Omega$  (cf. \cite{EchartLancaster1})  and $\nu$  is the exterior unit normal to $\partial\Omega.$

For a smooth domain, some type of boundary curvature condition (which depends on $H$)  must be satisfied in order to guarantee that a classical solution 
of (\ref{eq:D})-(\ref{bc:D}) exists for each $\phi\in C^{0}\left( \partial \Omega\right);$  when $H\equiv 0,$  this curvature condition is that $\partial\Omega$  
must have nonnegative curvature (with respect to the interior normal direction of $\Omega$) at each point (e.g. \cite{JenkinsSerrin}).
However, Leon Simon (\cite{Simon}) has shown that if $\Gamma_{0}\subset \partial\Omega$  is smooth (i.e. $C^{4}$), $H\equiv 0,$
$\phi\in C^{0,1}(\partial\Omega),$  the curvature $\Lambda$  of $\partial\Omega$  is negative on $\Gamma_{0}$  and $\Gamma$  is a compact subset of $\Gamma_{0},$
then the variational solution $z=f({\bf x}),$  ${\bf x}\in\Omega,$  extends to $\Omega\cup \Gamma$  as a H\"older continuous function with Lipschitz continuous trace, 
even though $f$  may not equal $\phi$  on $\Gamma;$  Simon's result holds for least area hypersurfaces in $\Real^{n},$  $n\ge 2$  when the mean curvature of 
$\partial\Omega$  has a negative upper bound on $\Gamma\subset \partial\Omega$  (see also \cite{Bour, Lin}).

One can look at this in a different way.  In the case $H\equiv 0,$  the requirement that $\Lambda({\bf p})<0$  at a point ${\bf p}\in\partial\Omega$  
implies that $Nf=0$  has a (continuous) Bernstein function $\psi$  at ${\bf p}$  for $\Omega$  (see Definition (\ref{Bernstein+})  and 
Definition (\ref{Bernstein-})).    
In \cite{EL1986B}, Bernstein functions for the minimal surface equation in $\Real^{2}$  are constructed for $C^{2,\alpha}$  domains $\Omega\subset\Real^{2}$  
whose curvature $\Lambda$  (with respect to $-\nu$)   vanishes at a finite number of points and satisfies $\Lambda\le 0$  on a segment of $\partial\Omega.$   
Using these Bernstein functions, we will prove the following generalization of \cite{Simon} when $n=2.$  
\begin{cor}
\label{Cor1}
Let $\Omega$  be a  domain in $\Real^{2},$  $\Gamma$  is a $C^{2,\lambda}$  open subset of $\partial\Omega$ and the  curvature $\Lambda$  (with respect to $-\nu$)  
of $\Gamma$  is nonpositive and vanishes at only a finite number of points of $\Gamma,$   for some $\lambda\in (0,1).$  
Suppose $\phi\in L^{\infty}(\partial\Omega),$  ${\bf y}\in \Gamma,$  
either $f$  is symmetric with respect to a line through ${\bf y}$  or $\phi$  is continuous at ${\bf y},$  and $f\in BV(\Omega)$  minimizes 
\begin{equation}
\label{JJ}
J(u)=\int_{\Omega} \sqrt{1+|Du|^{2}} d{\bf x} + \int_{\partial\Omega} |u-\phi| ds
\end{equation}
for $u\in BV(\Omega).$  Then $f\in C^{0}(\Omega\cup \{ {\bf y}\}).$  
If $\phi\in C^{0}(\Gamma),$  then $f\in C^{0}(\Omega\cup\Gamma).$  
\end{cor} 

\begin{figure}[htb]
\centerline{
\includegraphics[width=1.5in]{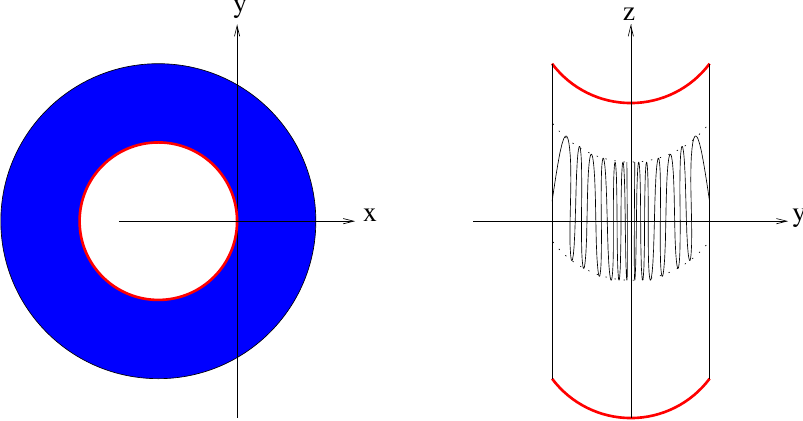}}
\caption{$\Omega$  and part of the graph of $\phi$  \label{Cat}}
\end{figure}

\begin{example}
\label{Example1}
Let $\Omega=\{(x,y)\in \Real^{2} : 1< (x+1)^{2}+y^{2}<\cosh^{2}(1)\}$  and $\phi(x,y)=\sin\left(\frac{\pi}{x^2+y^2}\right)$  for $(x,y)\neq (0,0)$  
(see Figure \ref{Cat} for a rough illustration of the graph of $\phi$).  Set ${\cal O}=(0,0).$
Let $f\in C^{2}(\Omega)$  minimize (\ref{JJ}) over $BV(\Omega)$  (i.e. $f$  is the variational solution of (\ref{eq:D})-(\ref{bc:D}) with $H \equiv 0$). 
Then Corollary \ref{Cor1} (with ${\bf y}={\cal O}$) implies $f\in C^{0}\left(\overline{\Omega}\right),$  even though $\phi$  has no limit at ${\cal O}.$  
\end{example}

Variational solutions of (\ref{eq:C})-(\ref{bc:C}) will exist in some sense (e.g. \S7.3 of \cite{FinnBook}) but they need not be finitely valued 
(e.g. the discussion of extremal curves in Chapter 6 of \cite{FinnBook}), bounded (e.g. \cite{FinnBook}, Corollary 5.5) or continuous at each point of the boundary 
(e.g. \cite{Korevaar}).  Variational solutions of (\ref{eq:D})-(\ref{bc:D}) will be bounded if $\phi\in L^{\infty}(\Omega)$  but need not be continuous 
at each point of the boundary.  
Many authors (e.g. \cite{EL1986A,NoraKirk1,Ger2,Lan:89,NitscheBook,Simon,Williams}) have investigated the boundary behavior at 
corners of variational solutions of (\ref{eq:D})-(\ref{bc:D}) and a number of authors 
(e.g. \cite{CF:94,NoraKirk2,FinnBook,Finn:96,Giu178,Korevaar,CFC,NCFC,LS1,Shi}) 
have done so for variational solutions of (\ref{eq:C})-(\ref{bc:C}).  

We shall investigate the existence and behavior of the radial limits of nonparametric prescribed mean curvature surfaces at corners of the domain, 
including ``smooth corners'' (e.g. Corollary \ref{Cor1}).  In particular, we shall use Bernstein functions to investigate the behavior of variational solutions 
of  (\ref{eq:D})-(\ref{bc:D}) or  (\ref{eq:C})-(\ref{bc:C}) at points of $\partial\Omega.$  

\section{Radial Limit Theorems}

Let  $Q$  be the operator on $C^{2}(\Omega)$  given by 
\begin{equation}
\label{Q}
Qf({\bf x}) \myeq Nf({\bf x}) - 2H({\bf x},f({\bf x})), \ \ \ \ {\bf x}\in\Omega, 
\end{equation}
where $H:\Omega\times\Real\to\Real$  is prescribed and $H({\bf x},t)$  is weakly increasing in $t$  for each ${\bf x}\in\Omega.$
Let $\nu$  be the exterior unit normal to $\partial\Omega,$  defined almost everywhere on $\partial\Omega.$    
We assume that for almost every ${\bf y}\in\partial\Omega,$  there is a continuous extension $\hat\nu$  of $\nu$  to a neighborhood of ${\bf y}.$ 

For each point ${\bf y}\in\partial\Omega,$  polar coordinates relative to ${\bf y}$  are denoted by $r_{\bf y}$ and $\theta_{\bf y}.$ 
We shall assume that for each ${\bf y}\in\partial\Omega,$  there exists a $\delta>0$  such that 
$\partial \Omega \cap B_{\delta}({\bf y})\setminus \{{\bf y}\}$ consists of two (open) arcs  ${\partial}_{\bf y}^{1}\Omega$  and $\partial_{\bf y}^{2}\Omega,$  
whose tangent rays approach the rays $L^{1}_{\bf y}: \:  \theta_{\bf y} = \alpha({\bf y})$  and $L^{2}_{\bf y}: \:  \theta_{\bf y} = \beta({\bf y})$ 
respectively, as the point ${\bf y}$ is approached, with $\alpha({\bf y})<\beta({\bf y})<\alpha({\bf y})+2\pi,$  in the sense that the tangent cone 
to $\overline{\Omega}$  at ${\bf y}$  is $\{\alpha({\bf y})\le\theta_{\bf y} \le\beta({\bf y}), 0\le r_{\bf y}<\infty\}.$
(In particular, $\{\alpha({\bf y})<\theta_{\bf y} <\beta({\bf y}), 0<r_{\bf y}<\epsilon(\theta_{\bf y})\}$  
is a subset of $\Omega$  for some $\epsilon\in C^{0}((\alpha({\bf y}),\beta({\bf y}))),$  $\epsilon(\cdot)>0,$  
and $\{\beta({\bf y})<\theta_{\bf y} <\alpha({\bf y})+2\pi, 0<r_{\bf y}<\epsilon(\theta_{\bf y})\}\cap\Omega=\emptyset$  
for some $\epsilon\in C^{0}((\beta({\bf y}),\alpha({\bf y})+2\pi)),$  $\epsilon(\cdot)>0.$)
When $\beta({\bf y})-\alpha({\bf y}) < \pi,$   $\partial\Omega$  is said to have a {\bf convex corner} at ${\bf y}$ and when $\beta({\bf y})-\alpha({\bf y}) > \pi,$   
$\partial\Omega$  is said to have a {\bf nonconvex corner} at ${\bf y}.$  
The radial limit of $f$  at ${\bf y}=(y_{1},y_{2})\in\partial\Omega$  in the direction $\omega(\theta)=(\cos\theta,\sin\theta),$  
$\theta\in \left(\alpha({\bf y}),\beta({\bf y})\right),$   is 
\begin{equation}
\label{radial}
Rf(\theta,{\bf y}) \myeq \lim_{r\downarrow 0} f(y_{1}+r\cos(\theta),y_{2}+r\sin(\theta)). 
\end{equation}
$Rf(\alpha({\bf y}),{\bf y})$  will be defined as the limit at ${\bf y}$  of the trace of $f$  restricted to 
${\partial}_{\bf y}^{1}\Omega$  and $Rf(\beta({\bf y}),{\bf y})$  as the limit at ${\bf y}$  of the trace of $f$  restricted to 
${\partial}_{\bf y}^{2}\Omega.$  
Notice that if $f$  is a generalized (e.g. variational or Perron) solution of (\ref{eq:D})-(\ref{bc:D}), $f$  need not equal $\phi$  on portions of 
$\partial\Omega$  and the tangential radial limits $Rf(\alpha({\bf y}),{\bf y})$  and $Rf(\beta({\bf y}),{\bf y})$  may, for example, differ from
$\phi({\bf y})$  when $\phi$  is continuous at ${\bf y}.$  

\begin{figure}[ht]
\centering
\includegraphics{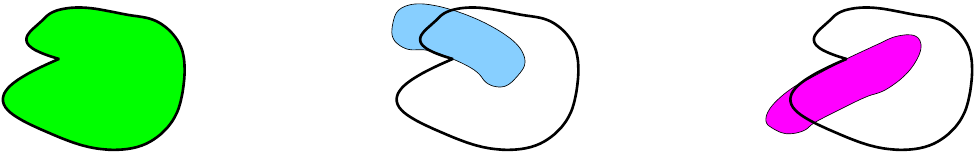}
\caption{$\Omega$ (left) \hspace{10mm} $U^{\pm}_{2}$ (middle)  \hspace{10mm} $U^{\pm}_{1}$ (right) \label{ONEa}}
\end{figure}

\begin{definition}
\label{Bernstein+} 
Given a domain $\Omega$  as above, a \underline {\bf upper Bernstein pair} $\left(U^{+},\psi^{+}\right)$  for a curve $\Gamma\subset \partial\Omega$  
and a function $H$  is 
a $C^{1}$  domain $U^{+}$  and a function $\psi^{+}\in C^{2}(U^{+})\cap C^{0}\left(\overline{U^{+}}\right)$  such that $\Gamma\subset\partial U^{+},$  
$\nu$  is the exterior unit normal 
to $\partial U^{+}$  at each point of $\Gamma$  (i.e. $U^{+}$  and $\Omega$  lie on the same side of $\Gamma;$  see Figure \ref{ONEa}), 
$Q\psi^{+}\le 0$  in $U^{+},$  and $T\psi^{+}\cdot\nu=1$  almost everywhere on an open subset of $\partial U^{+}$  containing $\overline{\Gamma}$  
in the same sense as in  \cite{CF:74a}; that is, for almost every ${\bf y}\in\Gamma,$ 
\begin{equation}
\label{zero}
\lim_{U^{+}\ni {\bf x}\to {\bf y}} \frac{\nabla \psi^{+}({\bf x})\cdot \hat\nu({\bf x})}{\sqrt{1+|\nabla \psi^{+}({\bf x})|^{2}}} = 1.
\end{equation} 
\end{definition}

\begin{definition}
\label{Bernstein-} 
Given a domain $\Omega$  as above, a \underline {\bf lower Bernstein pair} $\left(U^{-},\psi^{-}\right)$  for a curve $\Gamma\subset \partial\Omega$  
and a function $H$  is 
a $C^{1}$  domain $U^{-}$  and a function $\psi^{-}\in C^{2}(U^{-})\cap C^{0}\left(\overline{U^{-}}\right)$  such that $\Gamma\subset\partial U^{-},$  
$\nu$  is the exterior unit normal 
to $\partial U^{-}$  at each point of $\Gamma$  (i.e. $U^{-}$  and $\Omega$  lie on the same side of $\Gamma$), 
$Q\psi^{-}\ge 0$  in $U^{-},$  and $T\psi^{-}\cdot\nu=-1$  almost everywhere on an open subset of $\partial U^{-}$  containing $\overline{\Gamma}$  
in the same sense as in  \cite{CF:74a}.
\end{definition}

\noindent In the following theorem, we consider a domain with a nonconvex corner ${\bf y}$  and prove that the radial limits of $f$  at ${\bf y}$  
exist and behave as in \cite{EL1986A,Lan1985,Lan1988,LS1}.  In \cite{Lan1985}, $\Omega$  was required to be 
locally convex at points of ${\partial}_{\bf y}^{1}\Omega$  and $\partial_{\bf y}^{2}\Omega$  and, in \cite{EL1986A,Lan1988}, the curvatures of 
${\partial}_{\bf y}^{1}\Omega$  and $\partial_{\bf y}^{2}\Omega$   were required to have an appropriate positive lower bound when these curves were 
smooth.  In \cite{NoraKirk1}, no such curvature requirement was imposed but only nontangential radial limits were shown to exist.  
This theorem strengthens Theorem 1 of \cite{NoraKirk1} when the curvatures of ${\partial}_{\bf y}^{1}\Omega$  and $\partial_{\bf y}^{2}\Omega$ 
imply Bernstein functions exist (see \S \ref{BF}). 

\begin{thm}
\label{ONE}
Let $f\in C^{2}(\Omega)\cap L^{\infty}(\Omega)$  satisfy $Qf=0$  in $\Omega$  and let $H^{*}\in L^{\infty}(\Real^{2})$  satisfy  
$H^{*}({\bf x})=H({\bf x},f({\bf x}))$  for ${\bf x}\in\Omega.$
Suppose that ${\bf y}\in\partial\Omega,$  $\beta({\bf y})-\alpha({\bf y}) > \pi,$  
and there exist $\delta>0$  and upper and lower Bernstein pairs $\left(U^{\pm}_{1},\psi^{\pm}_{1}\right)$  and $\left(U^{\pm}_{2},\psi^{\pm}_{2}\right)$  for 
$(\Gamma_{1},H^{*})$  and $(\Gamma_{2},H^{*})$  respectively, where $\Gamma_{1}=B_{\delta}({\bf y})\cap {\partial}_{\bf y}^{1}\Omega$  and 
$\Gamma_{2}=B_{\delta}({\bf y})\cap {\partial}_{\bf y}^{2}\Omega.$      
Then the limits 
\begin{equation}
\label{Apple}
\lim_{\Gamma_{1}\ni {\bf x}\to {\bf y}} f({\bf x})=z_{1} \ {\rm and} \ 
\lim_{\Gamma_{2}\ni {\bf x}\to {\bf y}} f({\bf x})=z_{2}
\end{equation}
exist,  the radial limit $Rf(\theta,{\bf y})$  
exists for each $\theta\in [\alpha({\bf y}),\beta({\bf y})],$  $Rf(\alpha({\bf y}),{\bf y})=z_{1},$  $Rf(\beta({\bf y}),{\bf y})=z_{2},$  
and $Rf(\cdot,{\bf y})$  is a continuous function on $[\alpha({\bf y}),\beta({\bf y})]$  which behaves in one of the following ways:

\noindent (i) $Rf(\cdot,{\bf y})=z_{1}$  is a constant function  and $f$  is continuous at ${\bf y}.$  

\noindent (ii) There exist $\alpha_{1}$ and $\alpha_{2}$ so that $\alpha({\bf y}) \leq \alpha_{1}
< \alpha_{2} \leq \beta({\bf y}),$  $Rf=z_{1}$  on $[\alpha({\bf y}), \alpha_{1}],$  $Rf=z_{2}$  on $[\alpha_{2}, \beta({\bf y})]$
and $Rf$ is strictly increasing (if $z_{1}<z_{2}$) or strictly decreasing  (if $z_{1}>z_{2}$)  on $[\alpha_{1}, \alpha_{2}].$    
 
\noindent (iii) There exist $\alpha_{1}, \alpha_{L}, \alpha_{R}, \alpha_{2}$ so that
$\alpha({\bf y}) \leq \alpha_{1} < \alpha_{L} < \alpha_{R} < \alpha_{2} \leq \beta({\bf y}),$  
$\alpha_{R}= \alpha_{L} + \pi$, and $Rf$ is constant on $[\alpha({\bf y}), \alpha_{1}],
[ \alpha_{L}, \alpha_{R}]$, and $[ \alpha_{2}, \beta({\bf y})]$ and either  strictly increasing
on $[\alpha_{1}, \alpha_{L}]$ and  strictly decreasing on $[ \alpha_{R}, \alpha_{2}]$ or
strictly decreasing on $[\alpha_{1}, \alpha_{L}]$ and  strictly increasing on $[\alpha_{R},\alpha_{2}]$.  
\end{thm}

\noindent In the second theorem, we consider a domain with a smooth corner ${\bf y}$  (i.e. $\beta({\bf y})-\alpha({\bf y}) = \pi$)  
and show that the radial limits of $f$  at ${\bf y}$  exist and behave as expected.  
Corollary \ref{Cor1} follows from this theorem and an additional argument.

\begin{thm}
\label{TWO}
Let $f\in C^{2}(\Omega)\cap L^{\infty}(\Omega)$  satisfy $Qf=0$  in $\Omega$  and let $H^{*}\in L^{\infty}(\Real^{2})$  satisfy  
$H^{*}({\bf x})=H({\bf x},f({\bf x}))$  for ${\bf x}\in\Omega.$
Suppose that ${\bf y}\in\partial\Omega,$  $\beta({\bf y})-\alpha({\bf y}) = \pi,$  
and there exist $\delta>0$  and upper and lower Bernstein pairs $\left(U^{\pm},\psi^{\pm}\right)$  for $(\Gamma,H^{*}),$  
where $\Gamma=B_{\delta}({\bf y})\cap {\partial}\Omega.$  
Then the limits 
\[
\lim_{\Gamma_{1}\ni {\bf x}\to {\bf y}} f({\bf x})=z_{1} \ {\rm and} \ 
\lim_{\Gamma_{2}\ni {\bf x}\to {\bf y}} f({\bf x})=z_{2}
\]
exist, $Rf(\theta,{\bf y})$  exists for each $\theta\in [\alpha({\bf y}),\beta({\bf y})],$
$Rf(\cdot,{\bf y})\in C^{0}([\alpha({\bf y}),\beta({\bf y})])$ , $Rf(\alpha({\bf y}),{\bf y})=z_{1},$  $Rf(\beta({\bf y}),{\bf y})=z_{2},$  
and $Rf(\cdot,{\bf y})$  behaves as in (i) or (ii) of Theorem \ref{ONE}.
\end{thm}

\noindent In the third theorem, we consider a domain with a convex corner ${\bf y}$  and prove that the radial limits of $f$  at ${\bf y}$  exist 
and behave as expected.  This theorem strengthens Theorem 2 of \cite{NoraKirk1}.

\begin{thm}
\label{THREE}
Let $f\in C^{2}(\Omega)\cap L^{\infty}(\Omega)$  satisfy $Qf=0$  in $\Omega$  and let $H^{*}\in L^{\infty}(\Real^{2})$  satisfy  
$H^{*}({\bf x})=H({\bf x},f({\bf x}))$  for ${\bf x}\in\Omega.$
Suppose that ${\bf y}\in\partial\Omega$  and there exist $\delta>0$  and upper and lower Bernstein pairs $\left(U^{\pm}_{2},\psi^{\pm}_{2}\right)$  
for $(\Gamma_{2},H^{*}),$  where $\Gamma_{2}=B_{\delta}({\bf y})\cap {\partial}_{\bf y}^{2}\Omega.$  
Suppose further that $z_{1}=\lim_{\Gamma_{1}\ni {\bf x}\to {\bf y} } f\left({\bf x}\right)$  exists, where 
$\Gamma_{1}=B_{\delta}({\bf y})\cap {\partial}_{\bf y}^{1}\Omega.$ 
Then 
\[ 
\lim_{\Gamma_{2}\ni {\bf x}\to {\bf y}} f({\bf x})=z_{2}
\]
exists, $Rf(\theta,{\bf y})$  exists for each $\theta\in [\alpha({\bf y}),\beta({\bf y})],$
$Rf(\cdot,{\bf y})\in C^{0}([\alpha({\bf y}),\beta({\bf y})])$ , $Rf(\alpha({\bf y}),{\bf y})=z_{1},$  $Rf(\beta({\bf y}),{\bf y})=z_{2},$  
and $Rf(\cdot,{\bf y})$  behaves as in (i), (ii) or (iii) of Theorem \ref{ONE}.
\end{thm}

\noindent In the fourth theorem, we generalize Theorem 2 of \cite{NoraKirk2}.

\begin{thm}
\label{FOUR}
Let $f\in C^{2}(\Omega)$  satisfy  $Qf=0$  in $\Omega$  and let $H^{*}\in L^{\infty}(\Real^{2})$  satisfy  
$H^{*}({\bf x})=H({\bf x},f({\bf x}))$  for ${\bf x}\in\Omega.$  
Suppose that ${\bf y}\in\partial\Omega,$  $\beta({\bf y})-\alpha({\bf y}) < \pi,$  
and there exist $\delta>0$  and upper and lower Bernstein pairs $\left(U^{\pm}_{2},\psi^{\pm}_{2}\right)$  for $(\Gamma_{2},H^{*}),$  where 
$\Gamma_{2}=B_{\delta}({\bf y})\cap {\partial}_{\bf y}^{2}\Omega.$  
Suppose further that $f\in C^{1}\left(\Omega\cup \partial^{1}_{{\bf y}}\Omega\cup \partial^{2}_{{\bf y}}\Omega\right),$   
$Tf({\bf x})\cdot \nu({\bf x})= \cos(\gamma({\bf x}))  \mbox{\ for \ } {\bf x}\in \partial^{1}_{{\bf y}}\Omega,$
and 
\[
\lim_{\partial^{1}_{{\bf y}}\Omega\ni {\bf x}\to {\cal O} } \gamma\left({\bf x}\right)=\gamma_{2}.
\]
Suppose also that there exist $\lambda_{1},\lambda_{2}\in [0,\pi]$  with $0<\lambda_{2}-\lambda_{1}<2\left(\beta({\bf y})-\alpha({\bf y})\right)$  
such that $\lambda_{1}\le \gamma({\bf x})\le \lambda_{2}$  
for ${\bf x}\in \partial^{2}_{{\bf y}}\Omega$  and $\pi-2\alpha-\lambda_{1}<\gamma_{2}<\pi+2\alpha-\lambda_{2}.$   
Then the conclusions of Theorem \ref{THREE} hold.
\end{thm}

\noindent In the fifth theorem, we generalize Theorem 1 of \cite{LS1} at the cost of extra boundary assumptions; Theorem 1 of \cite{CEL1} also 
generalizes the Lancaster-Siegel theorem but only obtains nontangential radial limits while here the existence of all radial limits is established 
while not requiring the contact angle to be bounded away from zero or $\pi.$  

\begin{thm}
\label{FIVE}
Let $f\in C^{2}(\Omega)\cap L^{\infty}(\Omega)$  satisfy (\ref{eq:D})   and (\ref{bc:D}) almost everywhere on $\partial\Omega.$  
Let $H^{*}\in L^{\infty}(\Real^{2})$  satisfy  $H^{*}({\bf x})=H({\bf x},f({\bf x}))$  for ${\bf x}\in\Omega.$
Let ${\bf y}\in\partial\Omega$  and suppose  there exist $\delta>0$  and upper and lower Bernstein pairs 
$\left(U_{1}^{\pm},\psi^{\pm}_{1}\right)$  and $\left(U^{\pm}_{2},\psi^{\pm}_{2}\right)$  for $(\Gamma_{1},H^{*})$  and $(\Gamma_{2},H^{*})$  
respectively, where 
$\Gamma_{1}=B_{\delta}({\bf y})\cap {\partial}_{\bf y}^{1}\Omega$  and $\Gamma_{2}=B_{\delta}({\bf y})\cap {\partial}_{\bf y}^{2}\Omega.$      
If $\beta({\bf y})-\alpha({\bf y}) \le \pi,$  assume there exist constants $\underline{\gamma}^{\, \pm},
\overline{\gamma}^{\, \pm}, 0 \le \underline{\gamma}^{\, \pm} \leq \overline{\gamma}^{\, \pm} \le \pi,$   satisfying
\[ 
\pi - (\beta({\bf y})-\alpha({\bf y})) < \underline{\gamma}^{+} + \underline{\gamma}^{-} 
\]
\[
\le \overline{\gamma}^{\, +} +  \overline{\gamma}^{\, -} < \; \pi + \beta({\bf y})-\alpha({\bf y})
\]
such that $\underline{\gamma}^{\pm}\leq \gamma^{\pm}(s) \leq \overline{\gamma}^{\,  \pm}$
for all $s\in (0,s_{0}),$ for some $s_{0}>0.$
Then the conclusions of Theorem \ref{ONE} hold.  
\end{thm}

\begin{figure}[htb]
\centerline{
\includegraphics[width=3in]{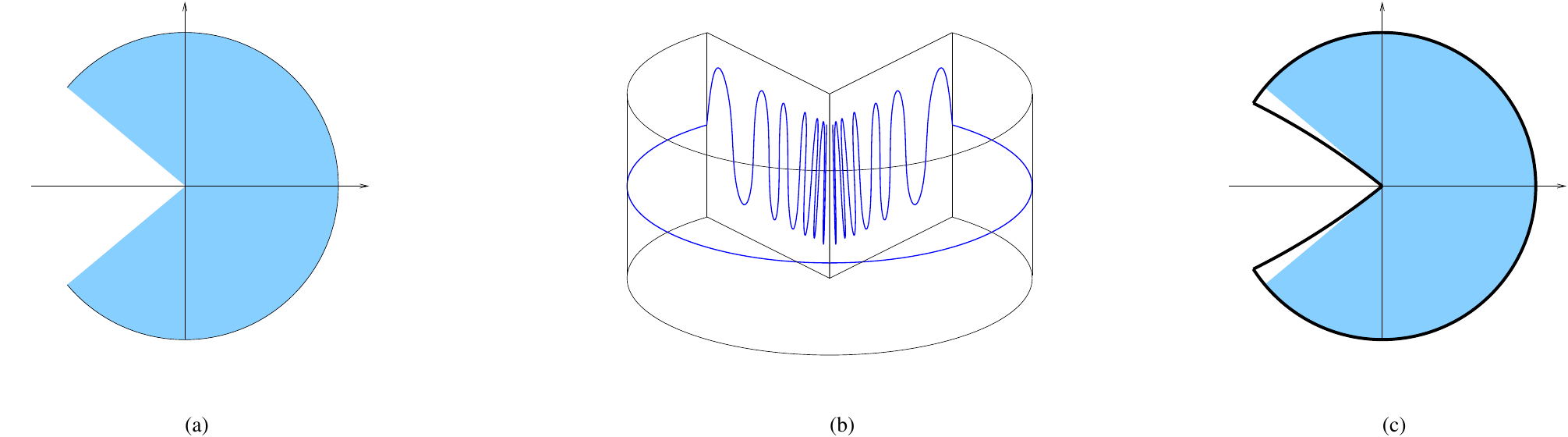}}
\caption{ (a) $\Omega$ \ \ \ (b) The graph of $\phi$  over $\partial\Omega$ \ \ \ (c) $\Omega_{1}$  \label{Dog}}
\end{figure}

\begin{example}
\label{Example2}
Let $\Omega = \{(r\cos\theta,r\sin\theta) : 0<r<1, -\alpha<\theta<\alpha \}$  with $\alpha>\frac{\pi}{2}.$     (see Figure \ref{Dog}(a)).  
Let  $\phi(x,y)=\sin\left(\frac{\pi}{x^2+y^2}\right)$  for $(x,y)\neq (0,0)$  (see Figure \ref{Dog}(b) for a rough illustration of the graph of $\phi.$).  
Let $f$  satisfy (\ref{eq:D})  in $\Omega$  with $H\equiv 0$  and $f=\phi$  on $\partial\Omega\setminus \{ {\cal O} \}.$  
Then \cite{NoraKirk1} shows that $Rf(\theta)$  exists when $|\theta|<\alpha.$  
Since $\Omega$  is locally convex at each point of $\partial\Omega\setminus \{ {\cal O}\},$  we see that $f\in C^{0}(\overline{\Omega}\setminus \{ {\cal O}\})$  
and $f=\phi$  on $\partial\Omega\setminus \{ {\cal O}\}.$
Since $\phi$ has no limit at ${\cal O},$  $Rf(\pm \alpha)$  do not exist; however $\lim_{\theta\downarrow -\alpha} Rf(\theta)$  and 
$\lim_{\theta\uparrow \alpha} Rf(\theta)$  both exist (e.g. from the behavior of $Rf(\theta)$  established in \cite{NoraKirk1,Lan1988,LS1}) 
and, by symmetry, are equal.  

Suppose we replace $\Omega$  with a slightly larger (and still symmetric) domain $\Omega_{1},$  $\Omega\subset\Omega_{1}\subset B_{1}({\cal O}),$
such that $\partial\Omega_{1}\cap B_{1}({\cal O})$  has negative curvature (with respect to the exterior normal to $\Omega_{1}$)  
and $\partial\Omega$  and $\partial\Omega_{1}$  are tangent at ${\cal O}$  (see Figure \ref{Dog} (c) for an illustration of $\Omega_{1}$).  
Let $f_{1}\in C^{2}(\Omega)$  minimize (\ref{JJ}) over $BV(\Omega_{1}),$  so that $f_{1}$  is the variational solution of (\ref{eq:D})-(\ref{bc:D})   
in $\Omega_{1}$  with $H\equiv 0.$    
Then Theorem \ref{ONE} implies $Rf_{1}(\theta)$  exists when $|\theta|\le \alpha$  and symmetry implies $Rf_{1}(-\alpha)=Rf_{1}(\alpha).$  
One wonders, for example, about the relationship between $Rf_{1}(\alpha)$  and $\lim_{\theta\uparrow \alpha} Rf(\theta).$
\end{example}

\section{Proofs}

\begin{rem}
The proofs of these Theorems are similar to those in \cite{NoraKirk1} (and \cite{CEL1}).  One difference is that the results in \cite{CEL1,NoraKirk1} were only 
concerned with nontangential radial limits at one point, ${\cal O},$  and so restricting the solution (``$f$'') to a subdomain which is tangent to the 
domain $\Omega$  at ${\cal O}$  and therefore assuming $f\in C^{0}(\overline{\Omega}\setminus \{ {\cal O} \})$  caused no difficulties.  
Since we wish to show that tangential radial limits also exist and describe the behavior of $f$  on $\partial\Omega,$  we cannot make such 
simplifying assumptions and so we have to modify the proofs in  \cite{CEL1,NoraKirk1}.
\end{rem}
\vspace{3mm}

\noindent {\bf Proof of Theorem \ref{ONE}:}  
We may assume $\Omega$  is a bounded domain.  Set  $S_{0} = \{ ({\bf x},f({\bf x})) : {\bf x} \in \Omega \}.$
From the calculation on page 170 of \cite{LS1},  we see that the area of $S_{0}$  is finite; let $M_{0}$  denote this area. 
For $\delta\in (0,1),$  set 
\[
p(\delta) = \sqrt{\frac{8\pi M_{0}}{\ln\left(\frac{1}{\delta}\right)}}.
\]
Let $E= \{ (u,v) : u^{2}+v^{2}<1 \}.$ 
As in \cite{EL1986A,LS1}, there is a parametric description of the surface $S_{0},$  
\begin{equation}
\label{PARAMETRIC}
Y(u,v) = (a(u,v),b(u,v),c(u,v)) \in C^{2}(E:{\Real}^{3}), 
\end{equation}
which has the following properties:
 
\noindent $\left(a_{1}\right)$  $Y$ is a diffeomorphism of $E$ onto $S_{0}$.
 
\noindent $\left(a_{2}\right)$
Set $G(u,v)=(a(u,v),b(u,v)),$  $(u,v)\in E.$   Then $G \in C^{0}(\overline{E} : {\Real}^{2}).$  

\noindent $\left(a_{3}\right)$  
%For each ${\bf y}\in\partial\Omega,$  
Set $\sigma({\bf y})=G^{-1}\left(\partial \Omega\setminus \{ {\bf y} \}\right);$  
then $\sigma({\bf y})$ is a connected arc of $\partial E$  and $Y$ maps $\sigma({\bf y})$  onto $\partial \Omega\setminus \{ {\bf y} \}.$  
We may assume the endpoints of $\sigma({\bf y})$  are ${\bf o}_{1}({\bf y})$  and ${\bf o}_{2}({\bf y}).$  
(Note that ${\bf o}_{1}({\bf y})$  and ${\bf o}_{2}({\bf y})$  are not assumed to be distinct.)
 
\noindent $\left(a_{4}\right)$
$Y$ is conformal on $E$: $Y_{u} \cdot Y_{v} = 0, Y_{u}\cdot Y_{u} = Y_{v}\cdot Y_{v}$
on $E$.
 
\noindent $\left(a_{5}\right)$
$\triangle Y := Y_{uu} + Y_{vv} = H^{*}\left(Y\right)Y_{u} \times Y_{v}$  on $E$.
\vspace{2mm} 

\noindent Notice that for each $C\in\Real,$   $Q(\psi^{+}_{j}+C)=Q(\psi^{+}_{j})\le 0$  on $\Omega\cap U^{+}_{j}$ 
and $Q(\psi^{-}_{j}+C)=Q(\psi^{-}_{j})\ge 0$  on $\Omega\cap U^{-}_{j},$
$j=1,2,$  and so 
\begin{equation}
\label{Barrier+}
N(\psi^{+}_{j}+C)({\bf x}) \le 2H({\bf x},f({\bf x}))=Nf({\bf x}) \ \ {\rm for} \ \ {\bf x}\in\Omega\cap U^{+}_{j}, \ \ j=1,2
\end{equation}
and 
\begin{equation}
\label{Barrier-}
N(\psi^{-}_{j}+C)({\bf x}) \ge 2H({\bf x},f({\bf x}))=Nf({\bf x}) \ \ {\rm for} \ \ {\bf x}\in\Omega\cap U^{-}_{j}, \ \ j=1,2.
\end{equation}
Let $q$  denote a modulus of continuity for $\psi^{\pm}_{1}$  and $\psi^{\pm}_{2}.$  

Let $\zeta({\bf y})=\partial E\setminus\sigma({\bf y});$  then $G(\zeta({\bf y}))=\{{\bf y}\}$  and ${\bf o}_{1}({\bf y})$  and ${\bf o}_{2}({\bf y})$  
are the endpoints of $\zeta({\bf y}).$  
There exists a $\delta_{1}>0$  such that if ${\bf w}\in E$  and ${\rm dist}\left({\bf w}, \zeta({\bf y})\right)\le 2\delta_{1},$  %$2\sqrt{\delta_{1}},$  
then $G({\bf w})\in \left(U^{+}_{1}\cup U^{+}_{2}\right)\cap \left(U^{-}_{1}\cup U^{-}_{2}\right).$  
Now $T\psi^{\pm}_{j}\cdot \nu=\pm 1$  (in the sense of \cite{CF:74a}) almost everywhere on an open subset $\Upsilon^{\pm}_{j}$  of $\partial U^{\pm}_{j}$  
which contains $\overline{\Gamma_{j}};$  there exists a $\delta_{2}>0$  such that 
$\left(\partial U^{\pm}_{j} \setminus \Upsilon^{\pm}_{j}\right) \cap \{{\bf x}\in \Real^{2} : |{\bf x}-{\bf y}|\le 2p(\delta_{2})\}=\emptyset.$  
Set $\delta^{*}=\min\{\delta_{1},\delta_{2}\}$  and 
\[
V^{*}= \{ {\bf w}\in E : {\rm dist}({\bf w},\zeta({\bf y}))<\delta^{*} \}.
\]
Notice if ${\bf w}\in V^{*},$  then $G({\bf w})\in U^{+}_{1}\cup U^{+}_{2}$  and $G({\bf w})\in U^{-}_{1}\cup U^{-}_{2}.$
\vspace{3mm} 

\noindent {\bf Claim:}   $Y$  is uniformly continuous on $V^{*}$  and so extends to a continuous function on $\overline{V^{*}}.$   
\vspace{3mm} 

\begin{figure}[ht]
\centering
\includegraphics{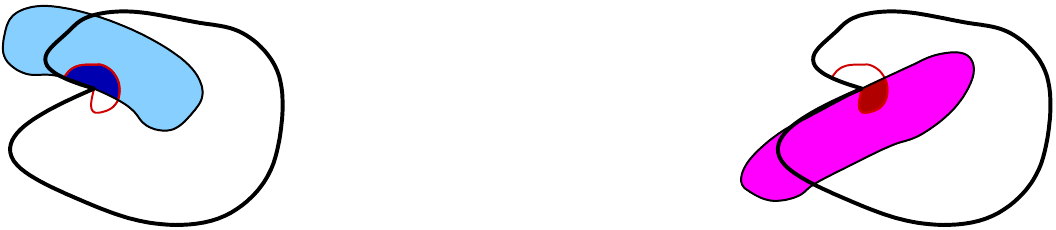}
\caption{General comparison principle applied on $U^{\pm}_{2}$ (left) and $U^{\pm}_{1}$ (right) \label{TWO_F}}
\end{figure}

\noindent {\bf Pf:}  Let $\epsilon>0.$  Choose $\delta\in \left(0,\left(\delta^{*}\right)^{2}\right)$  such that  $p(\delta)+2q(p(\delta))<\epsilon.$  
Let ${\bf w}_{1},{\bf w}_{2}\in V^{*}$  with $|{\bf w}_{1}-{\bf w}_{2}|<\delta;$  then 
$G({\bf w}_{1}), G({\bf w}_{2})\in \left(U^{+}_{1}\cup U^{+}_{2}\right) \cap \left(U^{-}_{1}\cup U^{-}_{2}\right) .$
Set $C_{r}({\bf w}) = \{ {\bf u} \in E : |{\bf u} - {\bf w}| = r \}$  and  $B_{r}({\bf w}) = \{ {\bf u} \in E : |{\bf u} - {\bf w}| < r \}.$
From the Courant-Lebesgue Lemma (e.g. Lemma $3.1$ in \cite{Cour:50}), we see that there exists $\rho=\rho(\delta)\in \left(\delta,\sqrt{\delta}\right)$  
such that the arclength $l_{\rho}({\bf w}_{1})$  of $Y(C_{\rho}({\bf w}_{1}))$  is less than $p(\delta).$  
Notice that ${\bf w}_{2}\in B_{\rho(\delta)}({\bf w}_{1}).$   Let 
\[
k(\delta)({\bf w}_{1})= \inf_{{\bf u}\in C_{\rho(\delta)}({\bf w}_{1})}c({\bf u}) = \inf_{ {\bf x}\in G(C_{\rho(\delta)}({\bf w}_{1})) } f({\bf x})
\]
and 
\[
m(\delta)({\bf w}_{1})= \sup_{{\bf u}\in C_{\rho(\delta)}({\bf w}_{1})}c({\bf u}) = \sup_{ {\bf x}\in G(C_{\rho(\delta)}({\bf w}_{1})) } f({\bf x});
\]
then  
\[
m(\delta)({\bf w}_{1})-k(\delta)({\bf w}_{1})\le l_{\rho} < p(\delta).
\]
Fix ${\bf x}_{0}\in C^{\prime}_{\rho(\delta)}({\bf w}_{1}).$  For $j=1,2,$  
set 
\[
C^{+}_{j}=\inf_{{\bf x}\in U^{+}_{j}\cap C^{\prime}_{\rho(\delta)}({\bf w}_{1}) } \psi^{+}_{j}({\bf x}) \ \ {\rm and} \ \ 
C^{-}_{j}=\sup_{{\bf x}\in U^{-}_{j}\cap C^{\prime}_{\rho(\delta)}({\bf w}_{1}) } \psi^{-}_{j}({\bf x}).
\]
Then $\psi^{+}_{j}-C^{+}_{j}\ge 0$  on $U^{+}_{j}\cap C^{\prime}_{\rho(\delta)}({\bf w}_{1})$  and 
$\psi^{-}_{j}-C^{-}_{j}\le 0$  on $U^{-}_{j}\cap C^{\prime}_{\rho(\delta)}({\bf w}_{1}).$ 
Therefore, for $j,l \in \{1,2\}$  and ${\bf x}\in U^{+}_{j}\cap U^{-}_{l}\cap C^{\prime}_{\rho(\delta)}({\bf w}_{1}),$  we have 
\[
k(\delta)({\bf w}_{1})+\left(\psi^{-}_{l}({\bf x})-C^{-}_{l}\right) \le f({\bf x}) \le m(\delta)({\bf w}_{1})+\left(\psi^{+}_{j}({\bf x})-C^{+}_{j}\right). 
\] 
For $j=1,2,$  set 
\[
b_{j}^{+}({\bf x})= m(\delta)({\bf w}_{1})+\left(\psi^{+}_{j}({\bf x})-C^{+}_{j}\right) \ \ \ \ 
{\rm for} \ \ {\bf x}\in U^{+}_{j}\cap \overline{G\left(B_{\rho(\delta)}({\bf w}_{1})\right)}
\]
and 
\[
b_{j}^{-}({\bf x})= k(\delta)({\bf w}_{1})+\left(\psi^{-}_{j}({\bf x})-C^{-}_{j}\right) \ \ \ \ 
{\rm for} \ \ {\bf x}\in  U^{-}_{j} \cap \overline{G\left(B_{\rho(\delta)}({\bf w}_{1})\right)}.
\]
Now $\rho(\delta)<\sqrt{\delta}<\delta^{*}\le\delta_{2};$  notice that if ${\bf w}\in \overline{B_{\rho(\delta)}({\bf w}_{1})},$  then 
$|{\bf w}-{\bf w}_{1}|<\delta_{2}$  and  $|G({\bf w})-{\bf y}|<2p(\delta_{2})$  and thus if 
${\bf x}\in \overline{G\left(B_{\rho(\delta)}({\bf w}_{1})\right)} \cap \partial U_{j}^{\pm},$  then ${\bf x}\in \Upsilon_{j}^{\pm}.$  
From (\ref{Barrier+}), (\ref{Barrier-}), the facts that $b_{l}^{-}\le f$  on $U^{-}_{l} \cap C^{\prime}_{\rho(\delta)}({\bf w}_{1})$  
and $f\le b_{j}^{+}$  on $U^{+}_{j} \cap C^{\prime}_{\rho(\delta)}({\bf w}_{1})$ 
for $j,l=1,2,$  and the general comparison principle (Theorem 5.1, \cite{FinnBook}), we have  (see Figure \ref{TWO_F}) 
\begin{equation}
\label{uniformA}
b_{l}^{-}\le f\ {\rm on} \ U^{-}_{l} \cap \overline{G\left(B_{\rho(\delta)}({\bf w}_{1})\right)} \ \ {\rm for} \ \  l=1,2
\end{equation}
and 
\begin{equation}
\label{uniformB}
f\le b_{j}^{+}\ {\rm on} \ U^{+}_{j}\cap \overline{G\left(B_{\rho(\delta)}({\bf w}_{1})\right)} \ \ {\rm for} \ \  j=1,2.
\end{equation}
\noindent Since the diameter of $G\left(B_{\rho(\delta)}({\bf w}_{1})\right)\le p(\delta),$  we have 
$\left|\psi^{\pm}_{j}({\bf x})-C^{\pm}_{j}\right|\le q(p(\delta))$  
for ${\bf x}\in U^{\pm}_{j}\cap G\left(B_{\rho(\delta)}({\bf w}_{1})\right).$  
Thus, whenever ${\bf x}_{1},{\bf x}_{2}\in \overline{G\left(B_{\rho(\delta)}({\bf w}_{1})\right)},$  at least one of the cases 
(a) ${\bf x}_{1},{\bf x}_{2} \in U^{+}_{1}\cap U^{-}_{1},$  (b)  ${\bf x}_{1},{\bf x}_{2} \in U^{+}_{2}\cap U^{-}_{2},$
(c) ${\bf x}_{1} \in U^{+}_{1}$  and ${\bf x}_{2}\in U^{-}_{2}$  or (d)  ${\bf x}_{1} \in U^{+}_{2}$  and ${\bf x}_{2}\in U^{-}_{1}$  holds.
Since $c({\bf w})=f\left(G({\bf w})\right),$  $G({\bf w}_{1})\in U^{+}_{i}\cap U^{-}_{j}$  for some $i=1,2$  and $j=1,2,$  and  
$G({\bf w}_{2})\in U^{+}_{l}\cap U^{-}_{n}$  for some $l=1,2$  and $n=1,2,$  we have 
\[
b_{j}^{-}\left(G({\bf w}_{1})\right)-b_{l}^{+}\left(G({\bf w}_{2})\right) \le c({\bf w}_{1})-c({\bf w}_{2}) \le 
b_{i}^{+}\left(G({\bf w}_{1})\right)-b_{n}^{-}\left(G({\bf w}_{2})\right)
\]
or 
\[
-\left[m(\delta)({\bf w}_{1})-k(\delta)({\bf w}_{1}) +\left(\psi^{+}_{l}(G({\bf w}_{2}))-C^{+}_{l}\right) - \left(\psi^{-}_{j}(G({\bf w}_{1}))+C^{-}_{j}\right)\right]
\]
\[
\le c({\bf w}_{1})-c({\bf w}_{2}) \le 
\]
\[
\left[m(\delta)({\bf w}_{1})-k(\delta)({\bf w}_{1}) +\left(\psi^{+}_{i}(G({\bf w}_{1}))-C^{+}_{i}\right) - \left(\psi^{-}_{n}(G({\bf w}_{2}))+C^{-}_{n}\right)\right].
\]
Since $|\psi^{\pm}_{j}(G({\bf w}))-C^{\pm}_{j}| \le q(p(\delta))$  for ${\bf w}\in B_{\rho(\delta)}({\bf w}_{1})\cap U^{\pm}_{j},$  we have   
\[
|c({\bf w}_{1})-c({\bf w}_{2})| \le  p(\delta)+2q(p(\delta))<\epsilon.
\]
Thus $c$  is uniformly continuous on  $V^{*}$  and, since $G \in C^{0}(\overline{E} : {\Real}^{2}),$  we see that  $Y$  is uniformly continuous on $V^{*}.$  
Therefore $Y$  extends to a continuous function, still denote $Y,$  on $\overline{V^{*}}.$ \qed
\vspace{3mm} 

Notice that 
\[
\lim_{\Gamma_{1}\ni {\bf x}\to {\bf y}} f({\bf x})=\lim_{\partial E\ni {\bf w}\to {\bf o}_{1}({\bf y})} c({\bf w})=c({\bf o}_{1}({\bf y}))
\]
and 
\[
\lim_{\Gamma_{2}\ni {\bf x}\to {\bf y}} f({\bf x})=\lim_{\partial E\ni {\bf w}\to {\bf o}_{2}({\bf y})} c({\bf w})=c({\bf o}_{2}({\bf y}))
\]
and so, with $z_{1}=c({\bf o}_{1}({\bf y}))$  and $z_{2}=c({\bf o}_{2}({\bf y})),$  we see that (\ref{Apple}) holds.

Now we  need to consider two cases: 
\begin{itemize}
\item[$\left(A\right)$] ${\bf o}_{1}({\bf y})= {\bf o}_{2}({\bf y}).$
\item[$\left(B\right)$] ${\bf o}_{1}({\bf y})\neq {\bf o}_{2}({\bf y}).$
\end{itemize}
These correspond to Cases 5 and 3 respectively in Step 1 of the proof of Theorem~1 of \cite{LS1}. 
\vspace{2mm} 

\noindent {\bf Case (A):}  Suppose ${\bf o}_{1}({\bf y})= {\bf o}_{2}({\bf y});$  set ${\bf o}={\bf o}_{1}({\bf y})= {\bf o}_{2}({\bf y}).$  
Then $f$  extends to a function in $C^{0}\left(\Omega \cup \{ {\bf y}\}\right)$  and case (i) of Theorem \ref{ONE} holds. 
\vspace{2mm} 

\noindent  {\bf Pf:} Notice that $G$  is a bijection of $E\cup \{{\bf o}\}$  and $\Omega\cup \{{\bf y}\}.$  Thus we may define $f=c\circ G^{-1},$  
so $f\left(G({\bf w})\right)=c({\bf w})$  for ${\bf w}\in E\cup \{{\bf o}\};$  this extends $f$  to a function defined on $\Omega\cup \{{\bf y}\}.$  
Let $\{\delta_{i}\}$  be a decreasing sequence of positive numbers converging to zero and consider the sequence of open sets 
$\{G(B_{\rho(i)}({\bf o}))\}$  in $\Omega,$  where $\rho(i)=\rho(\delta_{i}({\bf o})).$  
Now ${\bf y}\notin G(C_{\rho(i)}({\bf o}))$  and so there exist $\sigma_{i}>0$  such that 
\[
P(i)  =\{ {\bf x}\in\Omega : |{\bf x}-{\bf y}|<\sigma_{i}\} \subset G(B_{\rho(i)}({\bf o}))
\]
for each $i\in\Natural.$  Thus if  ${\bf x}\in P(i),$  we have $|f({\bf x}) - f({\bf y})|<p(\delta_{i})+2q(p(\delta_{i})).$  
The continuity of $f$  at ${\bf y}$  follows from this.  \qed 
\vspace{3mm} 

\noindent {\bf Case (B):}  Suppose ${\bf o}_{1}({\bf y})\neq {\bf o}_{2}({\bf y}).$   Then one of case (ii)  or (iii) of Theorem \ref{ONE} holds.
\vspace{2mm} 

\noindent  {\bf Pf:}
As at the end of Step 1 of the proof of Theorem 1 of \cite{LS1}, we define $X:B\to\Real^{3}$  by $X=Y\circ g$ and $K:B\to\Real^{2}$  
by $K=G\circ g,$  where $B=\{(u,v)\in\Real^{2} : u^{2}+v^{2}<1, \ v>0\}$  and $g:\overline{B}\to \overline{E}$  is either a conformal or an indirectly 
conformal (or anticonformal) map from $\overline{B}$  onto $\overline{E}$  such that  $g(1,0)= {\bf o}_{1}({\bf y}),$   
$g(-1,0)= {\bf o}_{2}({\bf y})$  and $g(u,0)\in  {\bf o}_{1}({\bf y}){\bf o}_{2}({\bf y})$  for each $u\in [-1,1],$  where 
${\bf ab}$  denotes the (appropriate) choice of arc in $\partial E$  with ${\bf a}$  and ${\bf b}$  as endpoints.  

Notice that $K(u,0)={\bf y}$  for $u\in [-1,1].$  Set $x=a\circ g,$  $y=b\circ g$  and $z=c\circ g,$  so that $X(u,v)=(x(u,v),y(u,v),z(u,v))$  
for $(u,v)\in B.$  Now, from Step 2 of the proof of Theorem 1 of \cite{LS1}, 
\[
X\in C^{0}\left(\overline{B}:\Real^{3}\right)\cap C^{1,\iota}\left(B\cup\{(u,0):-1<u<1\}:\Real^{3}\right)
\]
for some $\iota\in (0,1)$  and $X(u,0)=({\bf y},z(u,0))$  cannot be constant on any nondegenerate interval in $[-1,1].$  
Define $\Theta(u)= {\rm arg}\left( x_{v}(u,0)+iy_{v}(u,0) \right).$  From equation (12) of \cite{LS1}, we see that 
\[
\alpha_{1}=\lim_{u\downarrow -1} \Theta(u) \ \ \ \ {\rm and} \ \ \ \  \alpha_{2}=\lim_{u\uparrow 1} \Theta(u);
\]
here $\alpha_{1}<\alpha_{2}.$  
As in Steps 2-5 of the proof of Theorem 1 of \cite{LS1}, we see that $Rf(\theta)$  exists when $\theta\in \left(\alpha_{1},\alpha_{2}\right),$ 
\[
\overline{G^{-1}\left( L(\alpha_{2}) \right)} \cap \partial E = \{ {\bf o}_{1}({\bf y}) \} \    (\& \ 
\overline{K^{-1}\left( L(\alpha_{2}) \right)} \cap \partial B = \{ (1,0) \})   \  {\rm when}  \  \alpha_{2}<\beta({\bf y})
\]
\[
\overline{G^{-1}\left( L(\alpha_{1}) \right)} \cap \partial E = \{ {\bf o}_{2}({\bf y}) \} \ (\& \ 
\overline{K^{-1}\left( L(\alpha_{1}) \right)} \cap \partial B = \{ (-1,0) \})  \  {\rm when} \ \alpha_{1}>\alpha({\bf y})
\]
where  $L(\theta)= \{{\bf y}+(r\cos(\theta),r\sin(\theta))\in \Omega : 0<r<\delta^{*} \},$
and one of the following cases holds:
\vspace{1mm}

\noindent (a) $Rf$ is strictly increasing or strictly decreasing on $(\alpha_{1}, \alpha_{2})$.  
\vspace{1mm}

\noindent (b) There exist $\alpha_{L}, \alpha_{R}$ so that $\alpha_{1} < \alpha_{L} < \alpha_{R} < \alpha_{2},$  
$\alpha_{R}= \alpha_{L} + \pi$, and $Rf$ is constant on $[ \alpha_{L}, \alpha_{R}]$  and either increasing
on $(\alpha_{1}, \alpha_{L}]$ and decreasing on $[\alpha_{R}, \alpha_{2})$ or decreasing on $(\alpha_{1}, \alpha_{L}]$ 
and increasing on $[\alpha_{R}, \alpha_{2})$.  
\vspace{1mm}

\noindent We may argue as in Case A to see that $f$  is uniformly continuous on 
\[
\Omega^{+} =\{{\bf y}+(r\cos(\theta),r\sin(\theta))\in \Omega : 0<r<\delta, \alpha_{2}\le \theta<\beta({\bf y})+\epsilon\}
\]
and $f$  is uniformly continuous on 
\[
\Omega^{-} =\{{\bf y}+(r\cos(\theta),r\sin(\theta))\in \Omega : 0<r<\delta, \alpha({\bf y})-\epsilon< \theta\le\alpha_{1}\}
\]
for some small $\epsilon>0$  and $\delta>0,$  since  $G$  is a bijection of $E\cup \{{\bf o}_{1}({\bf y})\}$  and $\Omega\cup \{{\bf y}\}$  
and a bijection of $E\cup \{{\bf o}_{2}({\bf y})\}$  and $\Omega\cup \{{\bf y}\}.$  (Also see \cite{CEL1,NoraKirk2}.) 
Theorem \ref{ONE} then follows, as in \cite{NoraKirk1}, from Steps 2-5 of the proof of Theorem 1 of \cite{LS1} (replacing Step 3 with \cite{EchartLancaster1}).  
\qed
\vspace{3mm}

\noindent {\bf Proof of Theorem \ref{TWO}:}  The proof of this theorem is essentially the same as that of Theorem \ref{ONE}.
\vspace{3mm}

\noindent {\bf Proof of Corollary \ref{Cor1}:} From pp.1064-5 in \cite{EL1986B}, we see that there exist upper and lower Bernstein pairs 
$\left(U^{\pm},\psi^{\pm}\right)$  for $(\Gamma,H^{*}).$    
From Theorem \ref{TWO}, we see that the radial limits $Rf(\theta,{\bf y})$  exist for each $\theta\in [\alpha({\bf y}),\beta({\bf y})].$  
(Since $\beta({\bf y})-\alpha({\bf y})=\pi,$  case (iii) of Theorem \ref{ONE} cannot occur.)  
Set $z_{1}=Rf(\alpha({\bf y}),{\bf y}),$  $z_{2}=Rf(\beta({\bf y}),{\bf y})$  and $z_{3}=\phi({\bf y}).$  
If $z_{1}=z_{2},$  then case (i) of Theorem \ref{ONE} holds.  
(If $f$  is symmetric with respect to a line through ${\bf y},$  then $z_{1}=z_{2}$  and we are done.)

\begin{figure}[ht]
\centering
\includegraphics{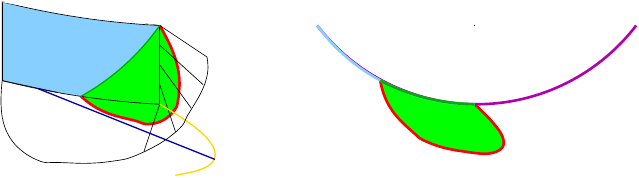}
\caption{Side View of $\Pi\cap \left(\Omega\times\Real\right)$ (left) and $\Omega_{1}$ (right)  \label{FOUR_F}}
\end{figure}

Suppose otherwise that $z_{1}\neq z_{2};$  we may assume that $z_{1}<z_{3}$  and $z_{1}<z_{2}.$    
Then there exist $\alpha_{1}, \alpha_{2}\in [\alpha({\bf y}),\beta({\bf y})]$  with $\alpha_{1}<\alpha_{2}$  such that 
\[
Rf(\theta,{\bf y}) \ \ {\rm is} 
\left\{ \begin{array}{ccc} 
{\rm constant}(=z_{1}) & {\rm for} & \alpha({\bf y})\le\theta\le\alpha_{1}\\ 
{\rm strictly \ increasing} & {\rm for} & \alpha_{1}\le\theta\le\alpha_{2}\\
{\rm constant}(=z_{2}) & {\rm for} & \alpha_{2}\le\theta\le\beta({\bf y}). \\ 
\end{array}
\right.
\]
From Theorem \ref{TWO}, we see that $Rf(\theta,{\bf y})$  exists for each $y\in \Gamma$  and $\theta\in [\alpha({\bf y}),\beta({\bf y})]$   
and $f$  is continuous on $\Omega\cup\Gamma\setminus \Upsilon$  for some countable subset $\Upsilon$  of $\Gamma.$  
Let $z_{0}\in \left(z_{1},\min\{z_{2},z_{3}\}\right)$  and $\theta_{0}\in (\alpha_{1},\alpha_{2})$  satisfy $Rf(\theta_{0},{\bf y})=z_{0}.$  
Let $C_{0}\subset\Omega$  be the $z_{0}-$level curve of $f$  which has ${\bf y}$  and a point ${\bf y}_{0}\in \partial\Omega\setminus \{{\bf y}\}$  
as endpoints. 
Let ${\bf y}_{1}\in \partial_{{\bf y}}^{1}\Omega\cap\Gamma\setminus \Upsilon$  and ${\bf y}_{2}\in C_{0}$  such that the (open) line segment $L$  joining
${\bf y}_{1}$  and ${\bf y}_{2}$  is entirely contained in $\Omega.$  
Let $M=\inf_{L} f,$  $\Pi$  be the plane containing $({\bf y},z_{0})$  and $L\times \{ M\},$  and $h$  be the affine function on $\Real^{2}$  whose graph is $\Pi.$  
Let $\Omega_{0}$  be the component of $\Omega\setminus \left(C_{0}\cup L\right)$  whose closure contains $B_{\delta}({\bf y})\cap \partial_{{\bf y}}^{1}\Omega$  
for some $\delta>0.$  Then there is a curve $C\subset \Omega_{0}$  on which $f=h$   whose endpoints are ${\bf y}_{3}$  and ${\bf y},$  for some 
${\bf y}_{3}\in \partial_{{\bf y}}^{1}\Omega$  between ${\bf y}_{1}$  and ${\bf y},$  such that $h>f$  in $\Omega_{1},$  where $\Omega_{1}\subset\Omega_{0}$  
is the open set bounded by $C$  and the portion of $\partial_{{\bf y}}^{1}\Omega$  between ${\bf y}$  and ${\bf y}_{3}.$    
Notice that $h<f$  in $L\cup C_{0}.$  
(In Figure \ref{FOUR_F}, on the left,  $\{\left({\bf x},h\left({\bf x}\right)\right): {\bf x}\in C\}$  is in red, $L$  is in dark blue, $C_{0}$  is in yellow, and 
the light blue region is a portion of $\partial_{{\bf y}}^{1}\Omega\times\Real,$  and, on the right, $\Omega_{0}$  is in light green and 
$\partial_{{\bf y}}^{2}\Omega$  is in magenta.) 
Now let $g\in C^{2}(\Omega)$  be defined by $g=f$  on $\Omega\setminus\overline{\Omega_{1}}$  and $g=h$  on $\Omega_{1}$  and observe that 
$J(g)<J(f),$  which contradicts the fact that $f$  minimizes $J.$  Thus it must be the case that $z_{1}=z_{2},$  case (i) of Theorem \ref{ONE} holds 
and $f$  is continuous at ${\bf y}.$  \qed  

\begin{rem}
\label{Remark8}
Corollary \ref{Cor1} can be generalized to minimizers of 
\[
J(u)=\int_{\Omega} \sqrt{1+|Du|^{2}} d{\bf x} + \int_{\Omega} \left( \int_{c}^{u({\bf x})} H({\bf x},t)  \ dt \right) d{\bf x}
+ \int_{\partial\Omega} |u-\phi| ds
\]
for $u\in BV(\Omega)$  and the conclusion remains the same; here $c$  is a reference height (e.g. $c=0$).  
In the proof of Corollary \ref{Cor1}, the only change is a replacement of the plane $\Pi$   with an appropriate surface (e.g. a portion of a sphere)  
over a subdomain like $\Omega_{1}$  such that the test function $g$  satisfies $J(g)<J(f).$ 
\end{rem}
\vspace{3mm} 

\noindent {\bf Proof of Example \ref{Example1}:} 
By Corollary \ref{Cor1}, $f$  is continuous on $\Omega\cup \{(0,0)\}.$  Clearly $f$  is continuous at $(x,y)$  when $(x+1)^2+y^2=\cosh^{2}(1).$ 
By \cite{Simon}, $f$  is continuous at $(x,y)$  when $(x+1)^2+y^2=1$  and $(x,y)\neq(0,0).$  
The parametrization (\ref{PARAMETRIC}) of the graph of $f$  (restricted to $\Omega\setminus \{(x,0):x<0\}$) satisfies $Y\in C^{0}(\overline{E}).$  
Notice that $\zeta((0,0))=\{ {\bf o} \}$  (since $\beta((0,0))-\alpha((0,0))=\pi$  and $z_{1}=z_{2}$)  for some ${\bf o}\in \partial E.$ 
Suppose $G$  in $(a_{2})$  is not one-to-one.  Then there exists a nondegenerate arc $\zeta\subset\partial E$  such that $G(\zeta)=\{{\bf y}_{1}\}$  for some 
${\bf y}_{1}\in \partial\Omega$  and therefore $f$  is not continuous at ${\bf y}_{1},$  which is a contradiction. 
Thus $f=g\circ G^{-1}$  and so $f\in C^{0}\left(\overline{\Omega}\right).$  
(The continuity of $G^{-1}$  follows, for example, from Lemma $3.1$ in \cite{Cour:50}.) 
\qed 
\vspace{3mm}

\noindent {\bf Proof of Theorem \ref{THREE}:}  The proof of Theorem 2 of \cite{NoraKirk1} uses unduloids as Bernstein functions (i.e. comparison surfaces) 
on subdomains of $\Omega$  (see Figure 7 of \cite{NoraKirk1}).  The proof of Theorem \ref{THREE} is essentially the same, using the 
Bernstein pairs $\left(U^{\pm},\psi^{\pm}\right)$ rather than unduloids, staying on $\partial_{\bf y}^{2}\Omega$  rather than on an arc of a circle inside 
$\Omega,$  and arguing as in the proof of Theorem \ref{ONE}.  \qed
\vspace{3mm} 

\noindent {\bf Proof of Theorem \ref{FOUR}:}  The proof of Theorem 2 of \cite{NoraKirk2} uses portions of tori as Bernstein functions (i.e. comparison surfaces) 
on subdomains of $\Omega$  (see Figure 7 of \cite{NoraKirk2}).  The proof of Theorem \ref{FOUR} is essentially the same, using the 
Bernstein pairs $\left(U^{\pm},\psi^{\pm}\right)$ rather than tori, staying on $\partial_{\bf y}^{2}\Omega$  rather than on an arc of a circle inside 
$\Omega,$  and arguing as in the proof of Theorem \ref{ONE}.  \qed
\vspace{3mm} 

\noindent {\bf Proof of Theorem \ref{FIVE}:}   The proof of Theorem 1 of \cite{CEL1} uses Theorem 2 of \cite{NoraKirk2}; the proof of Theorem \ref{FIVE} is 
essentially the same, using Theorem \ref{FOUR} in place of Theorem 2 of \cite{NoraKirk2} and arguing as in the proof of Theorem \ref{ONE}.  \qed
\vspace{3mm} 

\section{Bernstein Functions}
\label{BF}

The value of Theorems \ref{ONE} - \ref{FIVE} is dependent on the existence of Bernstein functions.  The results of \cite{EL1986B} %(pp.1063-5) 
provide Bernstein pairs for minimal surfaces.

\begin{prop}
\label{Prop1}
Let $a<b,$  $\lambda\in (0,1),$  $\psi\in C^{2,\lambda}([a,b])$  and  $\Gamma=\{(x,\psi(x))\in \Real^{2} : x\in [a,b]\}$   such that 
$\psi'(x)< 0$  for $x\in [a,b],$  $\psi''(x)<0$  for $x\in [a,b]\setminus J,$    there exist $C_{1}>0$  and $\epsilon_{1}>0$  such that 
if $\bar x\in J$  and $|x-\bar x|<\epsilon_{1},$  then $\psi''(x)\le -C_{1}|x-\bar x|^{\lambda},$  and 
$t\psi(x_{1})+(1-t)\psi(x_{2})<\psi\left(tx_{1}+(1-t)x_{2}\right)$  for each $t\in (0,1)$  and $x_{1},x_{2}\in [a,b]$  with $x_{1}\neq x_{2},$  
where $J$  is a finite subset of $(a,b).$
Then there exists an open set $U\subset\Real^{2}$  with  $\Gamma\subset \partial U$  and a function 
$h\in C^{2}\left(U\right)\cap C^{0}\left(\overline{U}\right)$  such that  $\partial U$  is a closed, $C^{2,\lambda}$  curve, 
$\Gamma$  lies below $U$  in $\Real^{2}$  (i.e. the exterior unit normal $\nu=(\nu_{1}(x),\nu_{2}(x))$  to $\partial U$ satisfies $\nu_{2}(x)<0$  
for $a\le x\le b$),   $Nh=0$  in $U$  and (\ref{zero}) holds for each ${\bf y}\in \Gamma,$  where $\hat\nu$ 
is a continuous extension  of $\nu$  to a neighborhood of $\Gamma.$
\end{prop}
\vspace{3mm} 

\noindent {\bf Proof:} We may assume that $a,b>0.$  
There exists $c>b$  and $k\in C^{2,\lambda}([-c,c])$  with $k(-x)=k(x)$  for $x\in [0,c]$  such that $k(x)=-\psi(x)$  for $x\in [a,b],$  
$k''(x)>0$  for $x\in [-c,c]\setminus J,$  where $J$  is a finite set, $k''(0)>0,$  and the set 
\[
K=\{(x,k(x))\in \Real^{2} : x\in [-c,c]\}  
\]
is strictly concave (i.e. $tk(x_{1})+(1-t)k(x_{2})>k\left(tx_{1}+(1-t)x_{2}\right)$  for each $t\in (0,1)$  and 
$x_{1},x_{2}\in [-c,c]$  with $x_{1}\neq x_{2}$).  
From \cite{EL1986B} (pp.1063-5), we can construct a domain $\Omega(K,l)$ such that $K\subset\partial\Omega(K,l)$  and 
$\Omega(K,l)$  lies below $K$  (i.e. the outward unit normal to  $\Omega(K,l)$  at $(x,k(x))$  is 
$\nu(x)=\frac{\left(-k'(x),1\right)}{\sqrt{1+\left(k'(x)\right)^{2}}};$  see Figure 4 of \cite{EL1986B}) 
and a function $F^{+}\in C^{2}\left(\Omega(K,l)\right)\cap C^{0}\left(\overline{\Omega(K,l)}\right)$  such that 
\[
\mu({\bf x})\myeq \frac{\left(\nabla F^{+}({\bf x}),-1\right)}{\sqrt{1+|\nabla F^{+}({\bf x})|^{2}}}, \ \ \ {\bf x}\in\Omega(K,l), 
\]
extends continuously to a function on $\Omega(K,l)\cup K$  and $\mu(x,k(x))\cdot \nu(x)=1$  for $x\in [-c,c].$  
Now let $V$  be an open subset of $\Omega$  with $C^{2,\lambda}$  boundary such that $\{(x,-\psi(x)) : x\in [a,b]\}\subset \partial V$  and 
$\partial V\cap \left(\partial\Omega(K,l)\setminus K\right)=\emptyset$  and then let 
$U=\{(x,-y) : (x,y)\in V\}$  and $h(x,y)=F^{+}(x,-y)$  for $(x,y)\in\overline{U}.$  \qed
\vspace{3mm}

\begin{rem}
\label{Remark2}
Let $\Omega\subset\Real^{2}$  be an open set,  $\Gamma\subset\partial\Omega$  be a $C^{2,\lambda}$  curve and ${\bf y}\in \Gamma$  
be a point at which we wish to have upper and lower Bernstein pairs for $H\equiv 0.$  
Let $\Sigma\subset\Gamma$  be the intersection of $\partial\Omega$  with a neighborood of ${\bf y}$  and suppose there is a rigid motion 
$\zeta:\Real^{2}\to\Real^{2}$  such that $\zeta\left(\Sigma\right)$  and $\zeta\left(\Omega\right)$  satisfy the hypotheses of Proposition \ref{Prop1}. 
Then $\left(\zeta^{-1}(U),h\circ\zeta\right)$  will be an upper Bernstein pair for $\Sigma$  and $H\equiv 0$  and 
$\left(\zeta^{-1}(U),-h\circ\zeta\right)$  will be a lower Bernstein pair for $\Sigma$  and $H\equiv 0.$  
\end{rem}
\vspace{3mm} 

When $H({\bf x},z)$  is independent of $z,$  the existence of (bounded) Bernstein functions is tied to boundary curvature conditions;   
in Theorem 3.1 of \cite{Giu178} (and Theorem 6.6 of \cite{FinnBook}), we see that 
\begin{prop}
\label{Prop2}
Suppose $\Omega$  is a $C^{2}$  domain in $\Real^{2}$  such that 
\begin{equation}
\label{curvature0}
|\int\int_{A} H({\bf x}) d{\bf x}| < \int |D\chi_{A}| \ \ \ {\rm for \ all\ } A\subset\Omega, \ A\neq \emptyset,\Omega
\end{equation}
and $\int\int_{\Omega} H({\bf x}) d{\bf x} = \int |D\chi_{\Omega}|;$  that is, $\Omega$  is an extremal domain.  
Let ${\bf y}\in\partial\Omega$  and suppose 
\begin{equation}
\label{curvature1}
\Lambda({\bf y})< 2H({\bf y}),
\end{equation}
where $\Lambda({\bf y})$  is the (signed) curvature of $\partial\Omega$  at ${\bf y}$  with respect to the interior normal direction.   
Then the (unique up to vertical translations) solution $g$  of $Ng({\bf x})=H({\bf x})$  for ${\bf x}\in\Omega$  is bounded and continuous in 
$W=\overline{\Omega}\cap B_{\epsilon}({\bf y}),$  $Tg$  extends continuously to a function on $W$  and $Tg({\bf x})=\nu({\bf x})$  for each 
${\bf x}\in B_{\epsilon}({\bf y}) \cap \partial\Omega$  for some $\epsilon>0,$    where $\nu$  is the exterior unit normal to $\Omega.$  
\end{prop}

\noindent Using Proposition \ref{Prop2} and a similar procedure to that in the proof of Proposition \ref{Prop1}, we can obtain Bernstein pairs near ${\bf y}$  
when $\partial\Omega\cap B_{\epsilon}({\bf y})$  is a subset of the boundary of an extremal domain $W$  for some $\epsilon>0$  such that $\Omega$  and $W$  
are on the same side of $\partial\Omega\cap B_{\epsilon}({\bf y})$   and the boundary curvature condition $\Lambda_{W}({\bf y})< 2|H({\bf y})|$  is satisfied.
In the same manner, we can obtain Bernstein pairs near ${\bf y},$  illustrated in Figure \ref{ONEa} by the sets $U^{\pm}_{1}$  and $U^{\pm}_{2},$ 
when ${\partial}_{\bf y}^{1}\Omega\cap B_{\epsilon}({\bf y})$  and ${\partial}_{\bf y}^{2}\Omega\cap B_{\epsilon}({\bf y})$  are subsets of the boundaries of 
extremal domains $W_{1}$  and $W_{2}$  for some $\epsilon>0,$   $\Omega$  and $W_{j}$  are on the same side of ${\partial}_{\bf y}^{j}\Omega\cap B_{\epsilon}({\bf y})$ 
for $j=1,2,$  $\Lambda_{W_{1}}({\bf y})< 2|H({\bf y})|$  and  $\Lambda_{W_{2}}({\bf y})< 2|H({\bf y})|,$ 
where $\Lambda_{W_{j}}({\bf y})$  denotes (signed) curvature of $\partial W_{j}.$  

\begin{rem}
\label{Remark3}
In Proposition \ref{Prop2},  the sets $A$  are Caccioppoli sets; that is, Borel sets such that the distributional (first) derivatives of the characteristic function 
$\chi_{A}$  of $A$  are Radon measures.  The notation $A\neq \emptyset,\Omega$  means that neither $A$  nor $\Omega\setminus A$  has (two-dimensional) measure zero and  
the notation $\int |D\chi_{\Omega}|$  means the total variation of $\chi_{A}\in BV(\Omega)$  (e.g. \S 6.3 of \cite{FinnBook}). 
Determining when hypothesis (\ref{curvature0}) is satisfied can be difficult; Giusti includes an Appendix in \cite{Giu178} which discusses the case of 
constant $H.$  
\end{rem}

We may use \S 14.4 of \cite{GT} (also see Corollary 14.13) to obtain Bernstein functions in a neighborhood $U$  of a point ${\bf y}\in\Gamma$  when  
$\Gamma\subset\partial\Omega$  is a $C^{2}$  curve satisfying $\Lambda({\bf x})<2|H({\bf x})|$  for ${\bf x}\in\Gamma\cap U$  and 
$H\in C^{0}\left(\overline{U\cap\Omega}\right)$  is either non-positive or non-negative in $U\cap\Omega.$  

\begin{lem}
\label{Lemma0}
Suppose $\Omega$  is a $C^{2,\lambda}$  domain in $\Real^{2}$  for some  $\lambda\in (0,1).$    Let ${\bf y}\in\partial\Omega$  and $\Lambda({\bf y})$  
denote the (signed) curvature of $\partial\Omega$  at ${\bf y}$  with respect to the interior normal direction (i.e. $-\nu$).  
Suppose $\Lambda({\bf y})< 2|H({\bf y})|$  and $H\in C^{0}\left(\overline{U\cap\Omega}\right)$  is either non-positive or non-negative in $U\cap\Omega,$  
where $U$  is some neighborhood of ${\bf y}.$  
Then there exist $\delta>0$  and upper and lower Bernstein pairs $\left(U^{\pm},\psi^{\pm}\right)$  for $(\Gamma,H),$  
where $\Gamma=B_{\delta}({\bf y})\cap {\partial}\Omega.$ 
\end{lem}
\vspace{3mm} 

\noindent {\bf Proof:} There exists $\delta_{1}>0$  such that $B_{\delta_{1}}({\bf y})\subset U$  and $\Lambda({\bf x})< 2|H({\bf x})|$  for each 
${\bf x}\in\partial\Omega \cap B_{\delta_{1}}({\bf y}).$   
There exists a $\delta_{2}\in (0,\delta_{1}/2)$  such that 
\[
\Lambda_{0}\myeq \sup\{\Lambda({\bf x}) : {\bf x}\in\partial\Omega\cap B_{\delta_{2}}({\bf y})\} < 
\inf\{2|H({\bf x})| : {\bf x}\in\partial\Omega\cap B_{\delta_{2}}({\bf y})\} \myeq 2H_{0}.
\]
If $\Lambda_{0}>0,$  set $R=\frac{1}{\Lambda_{0}};$  otherwise let $R$  be a small positive number.  
Now let $W$  be a $C^{2,\lambda}$  domain in $\Real^{2}$  such that $\partial\Omega\cap B_{\delta_{1}}({\bf y})\subset \partial W,$  
$\Omega$  and $W$  lie on the same side of $\partial\Omega\cap B_{\delta_{1}}({\bf y})$  and $W$  satisfies an interior sphere condition of radius $R$  
at each point of $\partial\Omega\cap B_{\delta_{2}}({\bf y}).$    
Continuously extend $H$  outside $U$  to $\overline{W}$  in such a manner that $H$  is either non-positive or non-negative in $W.$  
From inequality (14.73) of \cite{GT}, there exists $L>0$  such that 
\[
u({\bf x})-u_{0}({\bf x})\le L \ \ \ {\rm for} \ {\bf x}\in \partial W \cap B_{\delta_{2}}({\bf y}),
\]
where $u$  is any solution of (\ref{eq:D}) in $W$  and $u_{0}({\bf x})=\sup \{u({\bf t}) : {\bf t}\in\partial W\setminus B_{R}({\bf x}) \}.$  
We may assume $2\delta_{2}<R$  and set $u^{*}=\sup \{u({\bf t}) : {\bf t}\in\partial W\setminus B_{R-\delta_{2}}({\bf y}) \}.$
Then $u_{0}({\bf x})\le u^{*}$  for each ${\bf x}\in \partial W \cap B_{\delta_{2}}({\bf y})$  and 
$u\le L+u^{*}$  on $\partial\Omega\cap B_{\delta_{2}}({\bf y}).$  
Now let $\phi\in C^{\infty}(\partial W)$  such that $\phi=0$  on $\partial W\setminus B_{R-\delta_{2}}({\bf y})$  and $\phi>L$  on 
$\partial W \cap B_{\delta_{2}}({\bf y})$  and let $h\in C^{2}(W)$  be the solution of (\ref{eq:D})-(\ref{bc:D}) in $W$  with Dirichlet data $\phi.$  
(Just as \cite{GT} ignores in Theorem 14.11 the question of whether $u=\phi$  on $\partial\Omega\setminus B_{R}({\bf y}),$  we may assume that $W$  satisfies 
curvature conditions  (i.e. $\Lambda_{W}\ge 2|H|$)  on $\partial W\setminus B_{R-\delta_{2}}({\bf y})$  so that 
$h=\phi$  on $\partial\Omega\setminus B_{R}({\bf y})$  and so $h^{*}=0.$)   
It then follows (e.g. \cite{Bour}) that $h\in C^{0}(\overline{W})$  and 
\[
\frac{\partial h}{\partial\nu}=+\infty \ \ \ \ {\rm on} \ \ B_{\delta_{2}}({\bf y})\cap {\partial}\Omega.
\]
Thus $h$  is an upper Bernstein function.  The existence of a lower Bernstein function is similar.  \qed
\vspace{3mm} 

\begin{rem}
\label{Remark4}
In a similar manner, given ${\bf y}\in\partial\Omega$  we can establish the existence of upper and lower Bernstein pairs for the intersections 
of ${\partial}_{\bf y}^{1}\Omega$  and $\partial_{\bf y}^{2}\Omega$  with a neighborhood of ${\bf y}$  when these sets are each subsets of the 
boundaries of smooth (i.e. $C^{2,\lambda}$) domains $W_{1}$  and $W_{2}$  which satisfy appropriate boundary curvature conditions at ${\bf y}.$  
(For capillary surfaces in positive gravity (and prescribed mean curvature surfaces with $\frac{\partial H}{\partial z}({\bf x},z)\ge \kappa>0),$  
one can examine Theorem 2 of \cite{KorevaarSimon}.) 
\end{rem}

\section{Curvature Conditions on ${\partial}_{\bf y}^{1}\Omega$  and $\partial_{\bf y}^{2}\Omega$} 

In \cite{NoraKirk1}, the existence of nontangential radial limits of bounded, nonparametric prescribed mean curvature surfaces at nonconvex corners 
was proven; in Theorem \ref{ONE}, we showed that all radial limits of such surfaces at nonconvex corners exist when Bernstein functions exist.   
On the other hand, \cite{Lan:89} and Theorem 3 of \cite{LS1} provide examples in which no radial limit exists at a point ${\bf y}$  of $\partial\Omega$ 
at which the boundary of $\Omega$  is smooth.
In this section, we shall focus on the points ${\bf y}\in\partial\Omega$  at which $\beta({\bf y})-\alpha({\bf y}) \le \pi$  and ask which type of behavior 
(i.e. (a) no radial limits exist, (b) nontangential radial limits exist or (c) all radial limits exist) occurs, depending essentially on the curvatures of 
${\partial}_{\bf y}^{1}\Omega$  and $\partial_{\bf y}^{2}\Omega.$  
The following lemma shows that (a), (b) and (c) are the only possible behaviors of radial limits when $H({\bf x},t)$  is weakly increasing in $t$  
for each ${\bf x}\in\Omega,$  provided that we include in (b) all of the cases in which $Rf(\theta,{\bf y})$  exists for $\theta$  in one of the three intervals 
$(\alpha({\bf y}),\beta({\bf y})),$  $[\alpha({\bf y}),\beta({\bf y})),$  and $(\alpha({\bf y}),\beta({\bf y})].$ 

\begin{lem}
\label{Lemma1}
Let $f\in C^{2}(\Omega)\cap L^{\infty}(\Omega)$  satisfy $Qf=0$  in $\Omega$  and let $H^{*}\in L^{\infty}(\Real^{2})$  satisfy  
$H^{*}({\bf x})=H({\bf x},f({\bf x}))$  for ${\bf x}\in\Omega.$  
Let ${\bf y}\in\partial\Omega$  and suppose there exists a $\theta_{1}\in [\alpha({\bf y}),\beta({\bf y})]$  such that 
$Rf(\theta_{1},{\bf y})$  exists.  
Then $Rf(\theta,{\bf y})$  exists for each $\theta\in (\alpha({\bf y}),\beta({\bf y})),$      
$Rf(\cdot,{\bf y})\in C^{0}((\alpha({\bf y}),\beta({\bf y})))$  and $Rf(\cdot,{\bf y})$  behaves as in Theorem 1 of \cite{NoraKirk1}.  

Suppose, in addition, that there exist $\delta>0$  and upper and lower Bernstein pairs $\left(U^{\pm}_{1},\psi^{\pm}_{1}\right)$  and 
$\left(U^{\pm}_{2},\psi^{\pm}_{2}\right)$  for $(\Gamma_{1},H^{*})$  and $(\Gamma_{2},H^{*})$  respectively, where 
$\Gamma_{1}=B_{\delta}({\bf y})\cap {\partial}_{\bf y}^{1}\Omega$  and $\Gamma_{2}=B_{\delta}({\bf y})\cap {\partial}_{\bf y}^{2}\Omega.$      
Then the conclusions of Theorem \ref{ONE} hold.
\end{lem}
\vspace{3mm} 

\noindent {\bf Proof:}  The first part follows from Theorem 2 of \cite{NoraKirk1}.  The second part follows from Theorem \ref{THREE}. \qed
\vspace{3mm} 

Now suppose $f\in C^{2}(\Omega)$  satisfies  $Nf({\bf x})=H({\bf x})$  for ${\bf x}\in \Omega$  and ${\bf y}\in\partial\Omega$  satisfies 
$\beta({\bf y})-\alpha({\bf y}) \le \pi.$  Under what conditions do types of behavior (a), (b) or (c) occur?  
\vspace{3mm} 

\begin{lem}
\label{Lemma2}
Suppose  $\Omega,$  ${\bf y}$  and $H$  are as above and $\Lambda({\bf x})\ge 2|H({\bf x})|$  for almost all 
${\bf x}\in B_{\epsilon}({\bf y})\cap {\partial}_{\bf y}^{1}\Omega \cup \partial_{\bf y}^{2}\Omega,$  for some $\epsilon>0.$  
Then there exists $\phi\in L^{\infty}(\partial\Omega)$  such that the solution $f\in C^{2}(\Omega)\cap L^{\infty}(\Omega)$  of $Nf=H$  in $\Omega$  and 
$f=\phi$  almost everywhere on $\partial\Omega$  has no radial limits at ${\bf y}.$ 
\end{lem} 
\vspace{3mm} 

\noindent {\bf Proof:} This follows from Theorem 16.9 of \cite{GT} and the ``gliding hump'' argument in \cite{Lan:89}.  \qed  
\vspace{3mm} 

\begin{lem}
\label{Lemma3}
Suppose  $\Omega,$  ${\bf y}$  and $H$  are as above and $\Lambda({\bf x})< 2\left|H({\bf x})\right|$   for almost all 
${\bf x}\in B_{\epsilon}({\bf y})\cap {\partial}_{\bf y}^{1}\Omega \cup \partial_{\bf y}^{2}\Omega,$  for some $\epsilon>0.$ 
Then there exist $\delta>0$  and upper and lower Bernstein pairs 
$\left(U^{\pm}_{j},\psi^{\pm}_{j}\right)$  for $(\Gamma_{j},H),$  where $\Gamma_{j}=B_{\delta}({\bf y})\cap {\partial}_{\bf y}^{j}\Omega,$  for $j=1,2,$  
and the conclusions of Theorems \ref{TWO}--\ref{FIVE} hold when their other hypotheses are satisfied.  
\end{lem} 
\vspace{3mm} 

\noindent {\bf Proof:}  This follows from Remark \ref{Remark4}.  \qed
\vspace{3mm} 

\begin{thm}
\label{Theorem6}
Suppose $\Omega$  is a $C^{2,\lambda}$  domain in $\Real^{2}$  and $f\in C^{2}(\Omega)\cap L^{\infty}(\Omega)$  is a variational (i.e. BV) solution of  
(\ref{eq:D})-(\ref{bc:D}) for some $\phi\in L^{\infty}(\Omega)$  and $\lambda\in (0,1).$    
Let ${\bf y}\in\partial\Omega$  and let $\Lambda({\bf y})$  denote the (signed) curvature of $\partial\Omega$  at ${\bf y}$  with respect to the 
interior normal direction (i.e. $-\nu$).  
\begin{itemize}
\item[(i)]  Suppose $\Lambda({\bf y})< 2|H({\bf y})|.$  Then the conclusions of Theorem \ref{TWO} hold.
\item[(ii)]  Suppose $\Lambda({\bf y})> 2|H({\bf y})|.$  Then the conclusions of Theorem \ref{TWO} hold if $\phi$  restricted to ${\partial}_{\bf y}^{j}\Omega$ 
has a limit $z_{j}$  at ${\bf y}$  for $j=1,2,$  while for certain $\phi\in L^{\infty}(\Omega),$  $Rf(\cdot,{\bf y})$  does not exist 
for any $\theta\in [\alpha({\bf y}),\beta({\bf y})].$
\end{itemize}
\end{thm}
\vspace{3mm} 

\noindent {\bf Proof:}  The first part follows from Lemma \ref{Lemma3}. 
The second part follows from Theorem 16.9 of \cite{GT}, \cite{Lan1988} (see also \cite{EL1986A,LS1}) and Lemma \ref{Lemma2}.  \qed 
\vspace{3mm} 

\begin{rem}
\label{Remark5}
One can state a theorem similar to Theorem \ref{Theorem6} when $f\in C^{2}(\Omega)$  is a variational solution of (\ref{eq:D})-(\ref{bc:D}) 
for some $\phi\in L^{\infty}(\Omega)$  and ${\bf y}\in\partial\Omega$  satisfies $\beta({\bf y})-\alpha({\bf y})<\pi,$  ${\partial}_{\bf y}^{1}\Omega$  
and ${\partial}_{\bf y}^{2}\Omega$  are smooth, and boundary curvature conditions apply on ${\partial}_{\bf y}^{j}\Omega,$  $j=1,2.$   
If, for example, $\Lambda({\bf x})\ge 2|H({\bf x})|$  for ${\bf x}\in {\partial}_{\bf y}^{1}\Omega$  near ${\bf y},$  $\Lambda({\bf x})< 2|H({\bf x})|$  
for ${\bf x}\in {\partial}_{\bf y}^{2}\Omega$  near ${\bf y},$  and $\phi$  restricted to ${\partial}_{\bf y}^{1}\Omega$  has a limit $z_{1}$  at ${\bf y},$  
then the conclusions of Theorem \ref{TWO} hold  (e.g. Theorem 2 of \cite{NoraKirk1}).  
\end{rem}

\end{document}